\documentclass[letterpaper, 10 pt, conference]{ieeeconf}  

\usepackage[utf8]{inputenc}
\usepackage{amsfonts}
\usepackage{amsmath}
\usepackage{comment}
\usepackage{xcolor}
\usepackage{algorithm, algpseudocode}
\usepackage{subcaption}

\usepackage{pgfplots}
  \pgfplotsset{compat=newest}
  \usetikzlibrary{plotmarks}
  \usetikzlibrary{arrows.meta}
  \usepgfplotslibrary{patchplots}
  \usepackage{grffile}

\newcommand{\R}{\mathbb{R}}
\newcommand{\bmat}[1]{\begin{bmatrix}#1\end{bmatrix}}

\newcommand{\He}[1]{{\color{black}  #1}}

\newtheorem{prop}{Proposition}

\IEEEoverridecommandlockouts                              
\makeatletter
\let\NAT@parse\undefined
\makeatother
\usepackage{hyperref}
\title{\LARGE \bf
Continuous Abstraction of Nonlinear Systems using Sum-of-Squares Programming
}

\author{Stanley W. Smith, He Yin and Murat Arcak
\thanks{Funded in part by an NDSEG Graduate Fellowship, the Air Force Office of Scientific Research grant
FA9550-18-1-0253, and the Office of Naval Research grant N00014-18-1-2209.}
\thanks{The first two authors contributed equally to this work.}
\thanks{S. W. Smith and M. Arcak are with the Department of Electrical Engineering and Computer Sciences, University of California, Berkeley {\tt\small \{swsmth, arcak\}@eecs.berkeley.edu.}}%
\thanks{H. Yin is with the Department of Mechanical Engineering, University of California, Berkeley {\tt\small he\_yin@berkeley.edu.}}%
}

\begin{document}

\maketitle
\thispagestyle{empty}
\pagestyle{empty}

\begin{abstract}
We present a control design procedure for nonlinear control systems in which we represent a potentially high dimensional system with a low dimensional continuous-state abstraction. The abstraction generates a reference which the original system follows with a low level controller. We propose sum-of-squares programming as a tool to design this controller and to provide an upper bound on the relative error between the system and its abstraction. We compute the low level controller simultaneously with a simulation function that gives the boundedness guarantee for the relative error.
\end{abstract}

\section{Introduction}
An important task for cyber-physical systems is to design controllers that perform complex tasks while providing safety guarantees. Formal methods approaches \cite{belta2017formal} allow one to specify desired behaviors and then compute control strategies that achieve them, but are often limited to systems with small state dimension. Thus, it is beneficial to represent the system of interest with a lower dimensional model that is amenable to such approaches. A rigorous framework for such reduction is the notion of a continuous-state abstraction \cite{girard2009hierarchical}, henceforth called a \textit{continuous abstraction}, which provides a low level controller to ensure boundednesss of the error between the system and its abstraction as well as error bounds. In \cite{girard2009hierarchical}, the authors also provide conditions under which the error converges to zero, certified by a Lyapunov-like function referred to as a simulation function. This is extended in \cite{rungger2018compositional, zamani2018compositional} to interconnected control systems, where an additional condition on the interconnection topology is required, and the simulation function is constructed in a compositional manner using small-gain type reasoning and dissipativity properties of the subsystems, respectively.

In \cite{smith2018hierarchical}, an explicit connection is made between the aforementioned conditions and the existence of an invariant manifold on which the full system model reduces to the abstraction. In addition, \cite{smith2018hierarchical} and \cite{smith2018approximate} show that one can relax this invariance condition, resulting in a disturbance term in the dynamics relative to the manifold. This relaxation broadens the applicability of continuous abstraction; however, the boundedness analysis in \cite{smith2018hierarchical} and \cite{smith2018approximate} for the error system relies on an incremental dissipativity condition that may be hard to establish in practice.

In this paper, we do not require invariance or incremental dissipativity conditions; instead, we formulate the search for a low level controller and a simulation function for error bounds as a sum-of-squares (SOS) optimization problem. We approximate the system dynamics with a polynomial model and assume that the abstract state and input are restricted to semi-algebraic constraint sets. These sets and the model are then used as an input to our SOS optimization problem which, when feasible, returns a corresponding error bound. If the resulting error bound is too large, the constraint sets can be restricted further in an iterative procedure outlined in Section \ref{outline}.

Related publications other than those from the continuous abstraction literature include a modular control design approach to motion planning in \cite{herbert2017fastrack}, which consists of a planning layer and a tracking layer, and obtains an upper bound on the tracking error using Hamilton-Jacobi reachability analysis. An extension in \cite{singh2018robust} shows how this may also be done with SOS optimization. In \cite{kousik2017safe}, SOS optimization is used to construct a function that bounds the error between a high fidelity model and a low fidelity model. We note that these works focus on motion planning applications where the zero-error manifold is a subspace with a particular structure, whereas we allow any manifold to be used depending on the application. For example, in this paper we include a platooning application where some of the desired states are a translation of the abstract states, and in prior work \cite{smith2018approximate} we included a temperature control application where multiple states are aggregated into one abstract state. In \cite{Majumdar:17}, a bound is computed for the relative distance between the closed-loop system trajectory and a nominal trajectory that is given as a time series. However, a large number of bounds need to be precomputed for different nominal trajectories to consider all the possible maneuvers of the system.

The paper is organized as follows. In Section \ref{sec_overview}, the framework for continuous abstraction is reviewed, and methods for computing the error bounds and low level controllers are proposed. In Section \ref{example1} and \ref{example2}, our method is applied to two nonlinear system examples. Section \ref{sec_conclusion} summarizes the results.

\textit{Notation}: $\mathcal{C}^1$ is the set of differentiable functions whose derivative is continuous. For $\xi \in \mathbb{R}^n$, $\mathbb{R}[\xi]$ represents the set of polynomials in $\xi$ with real coefficients, and $\R^{m}[\xi]$ and $\R^{m\times p}[\xi]$ denote all vector and
matrix valued polynomial functions. The subset $\Sigma[\xi] := \{p = p_1^2 + p_2^2 + ... + p_M^2 : p_1, ..., p_M \in \mathbb{R}[\xi]\}$ of $\mathbb{R}[\xi]$ is the set of SOS polynomials in $\xi$. For $\gamma \in \mathbb{R}$, and continuous $V: \mathbb{R}^n \rightarrow \mathbb{R}$, we denote the sublevel set as
\begin{align}
\He{\Omega_{\gamma}^V := \{e \in \mathbb{R}^n : V(e) \leq \gamma \}.} \label{eq:sublevel}
\end{align}

\section{Continuous Abstraction} \label{sec_overview}
We consider control affine systems of the form
\begin{equation}
\Sigma: \quad \dot{x}(t) = f(x(t)) + g(x(t)) u(t), \label{concrete}
\end{equation}
where $x(t) \in \R^n, \ u(t) \in \R^m, \ f: \R^n \to \R^n$, and $g : \R^n \to \R^{n}\times \R^m$. We refer to the original system \eqref{concrete} as the concrete system. To facilitate the control design procedure, we represent \eqref{concrete} with a continuous abstraction
\begin{equation}
\hat{\Sigma}: \quad \dot{\hat{x}}(t) = \hat{f}(\hat{x}(t)) + \hat{g}(\hat{x}(t)) \hat{u}(t), \label{abstract}
\end{equation}
where $\hat{x}(t) \in \R^{\hat{n}}, \ \hat{u}(t) \in \R^{\hat{m}}, \ \hat{f} : \R^{\hat{n}} \to \R^{\hat{n}}$, and $\hat{g} : \R^{\hat{m}} \to \R^{\hat{n}} \times \R^{\hat{m}}$. Since the abstraction \eqref{abstract} typically has a smaller state dimension than the full system \eqref{concrete}, we define a map $\pi: \R^{\hat{n}} \to \R^n$ which allows comparison of $x(t)$ and $\hat{x}(t)$. In particular, we want $\Sigma$ and $\hat{\Sigma}$ to remain within close proximity to the manifold
\begin{equation} \label{manifold}
x(t) = \pi(\hat{x}(t))
\end{equation}
for all $t \geq 0$. Thus, we may interpret $\pi(\hat{x}(t))$ as generating a reference signal for $\Sigma$ to track. The examples in Sections \ref{example1} and \ref{example2} illustrated physically meaningful choices of $\pi$.

The abstract control $\hat{u}(t)$ could be designed formally so that \eqref{abstract} satisfies a high level objective expressible in, for example, Signal Temporal Logic (STL) \cite{maler2004monitoring}. Such objectives can be achieved with a numerical optimization-based approach, such as model predictive control (MPC) \cite{raman2017}. However, since such methods do not scale to large systems, a shrunken state dimension for \eqref{abstract} is beneficial.

\subsection{The Error System} \label{errorSystem}
To ensure boundedness of the error $e(t) := x(t) - \pi(\hat{x}(t))$, we incorporate a low level controller
\begin{equation} \label{interface}
u(t) = \kappa(e(t), \hat{x}(t), \hat{u}(t)) 
\end{equation}
to be designed and note that $e(t)$ evolves according to
\begin{align}
&\dot{e}(t) = f_e(e(t), \hat{x}(t), \hat{u}(t)) \label{errDyn} \\[2.5pt]
&\qquad + \He{g_e(e(t), \hat{x}(t))} \kappa(e(t), \hat{x}(t), \hat{u}(t)), \nonumber
\end{align}
\begin{align}
\text{where} \ &f_e(e(t), \hat{x}(t), \hat{u}(t)) := \label{f_e} \\
&\ \begin{aligned}
&f(\pi(\hat{x}(t)) + e(t)) - \frac{\partial \pi}{\partial \hat{x}} \cdot \left( \hat{f}(\hat{x}(t)) + \hat{g}(\hat{x}(t)) \hat{u}(t) \right),
\end{aligned} \nonumber \\[5pt]
&\He{g_e(e(t), \hat{x}(t))} := g(\pi(\hat{x}(t)) + e(t)). \label{g_e}
\end{align}
To design the low level controller (\ref{interface}) and to provide error bounds, we make use of Proposition \ref{prop1} below, which assumes that the abstract state and input are restricted to the sets $\hat{\mathcal{X}}$ and $\hat{\mathcal{U}}$, respectively. This assumption is satisfied by designing a control $\hat{u}$ which imposes constraints on $\hat{x}$ and $\hat{u}$.


\begin{prop} \label{prop1}
Given the error dynamics (\ref{errDyn}) with mappings $f_e: \mathbb{R}^n \times \mathbb{R}^{\hat{n}} \times \mathbb{R}^{\hat{m}} \rightarrow \mathbb{R}^n$, $g_e: \mathbb{R}^n \times \mathbb{R}^{\hat{n}} \rightarrow \mathbb{R}^n$, and $\gamma > 0$, $\hat{\mathcal{X}} \subseteq \mathbb{R}^{\hat{n}}$, $\hat{\mathcal{U}} \subseteq \mathbb{R}^{\hat{m}}$, if there exists a $\mathcal{C}^1$ function $V: \mathbb{R}^{n} \rightarrow \mathbb{R}^+$, and $\kappa: \mathbb{R}^n \times \mathbb{R}^{\hat{n}} \times \mathbb{R}^{\hat{m}} \rightarrow \mathbb{R}^m$, such that 
\begin{align}
&V(0) = 0 \ \text{and} \ V(e) > 0 \ \text{for all} \ e \neq 0,  \label{eq:A1} \\
&\frac{\partial V}{\partial e}\cdot(f_e(e, \hat{x}, \hat{u})+g_e(e, \hat{x})\kappa(e,\hat{x},\hat{u})) < 0, \nonumber \\
& \quad \quad \forall e, \hat{x}, \hat{u}, \ \text{s.t.} \ V(e) \ge \gamma, \ \hat{x} \in \hat{\mathcal{X}}, \ \hat{u} \in \hat{\mathcal{U}}, \label{eq:A2}
\end{align}
then the sublevel set $\Omega_{\gamma}^V$, defined in (\ref{eq:sublevel}), is forward invariant. If, further, $V$ is radially unbounded, then $\forall e(0) \in \R^n$, $e(t)$ remains bounded and converges to  $\Omega_{\gamma}^V$.
\end{prop}

\He{The set $\Omega_\gamma^V$ is an error bound for the error dynamics (\ref{errDyn}), associated with the control law $\kappa$.} If the low level control is subject to the constraints $u(t) \in \mathcal{U}$, where $\mathcal{U}:= \{u \in \R^m : \underline{u} \leq u \leq \overline{u}\}$, the symbol ``$\leq$" represents componentwise inequality, we augment the constraints (\ref{eq:A1}), (\ref{eq:A2}) and $V$ being radially  unbounded with the set containment constraints: for all $i = 1, ..., m$,
\begin{align}
    \Omega_\gamma^V \subseteq \{e \in \R^n : \kappa_i(e, \hat{x}, \hat{u}) \leq \overline{u}_i\}, \ \forall \hat{x} \in \hat{\mathcal{X}}, \ \hat{u} \in \hat{\mathcal{U}}, \label{eq:uUB}\\
    \Omega_\gamma^V \subseteq \{e \in \R^n : \underline{u}_i \leq \kappa_i(e, \hat{x}, \hat{u})\}, \ \forall \hat{x} \in \hat{\mathcal{X}}, \ \hat{u} \in \hat{\mathcal{U}}. \label{eq:uLB}
\end{align}

To achieve a less conservative bound, we formulate a high-level optimization problem with the volume of $\Omega_{\gamma}^V$ as the objective function and (\ref{eq:A1}) to (\ref{eq:uLB}) as constraints. We denote the problem as \textbf{opt}$\boldsymbol{(f_e, g_e, \gamma, \hat{\mathcal{X}}, \hat{\mathcal{U}}, \mathcal{U}})$:
	\begin{align}
	\min_{V, \kappa}\ &{\textnormal{volume}(\Omega_{\gamma}^V)}, \nonumber \\
	\textnormal{s.t.} \ &(\ref{eq:A1}) - (\ref{eq:uLB}), \text{and} \ V \ \text{is radially unbounded}. \nonumber
	\end{align}

\subsection{Sum-of-Squares formulation of the optimization problem} \label{SOSopt}
To formulate \textbf{opt} above as a SOS problem, assume that the sets $\hat{\mathcal{X}}$ and $\hat{\mathcal{U}}$ are semi-algebraic sets:
\begin{equation}
\begin{aligned}
    \hat{\mathcal{X}} :=& \{\hat{x} \in \R^{\hat{n}}: \hat{p}_1(\hat{x}) \leq 0 \},  \\
    \hat{\mathcal{U}} :=& \{\hat{u} \in \R^{\hat{m}}: \hat{p}_2(\hat{u}) \leq 0 \}, \label{stateInputConstraints}
\end{aligned}
\end{equation}
where the polynomial functions $\hat{p}_1$ and $\hat{p}_2$ are chosen based on the state and input constraints we enforced in the formal synthesis for the abstract system. We also assume that $f_e(e,\hat{x},\hat{u}) \in \R^n[(e,\hat{x},\hat{u})]$ and $g_e \in \R^{n\times m}[(e,\hat{x})]$. 

Using sum-of-squares relaxation for non-negativity, and the generalized S-procedure \cite{Parrilo:00} for set containment, we give the following SOS optimization problem denoted as \textbf{sosopt}$\boldsymbol{(f_e, g_e, \gamma, \hat{p}_1, \hat{p}_2, \underline{u}, \overline{u})}$: 
\begingroup
\allowdisplaybreaks
\begin{align}
    \min_{V, \kappa, s_i, s_k^j} \ & \textnormal{volume}(\Omega_{\gamma}^V) \nonumber \\
    \textnormal{s.t.} \ & s_i \in \Sigma[(e,\hat{x},\hat{u})], \ \forall i = 1,...,3,  \nonumber \\
    &s_k^j \in \Sigma[(e,\hat{x},\hat{u})], \ \forall j = 1,...,m, \ \forall k = 4,...,9, \nonumber \\
    & \kappa \in \R^m[(e, \hat{x}, \hat{u})], \ V(e)- L_1 \in \Sigma[e], \label{eq:A3} \\
    &-\frac{\partial V}{\partial e}\cdot(f_e + g_e \cdot \kappa) - L_2 + s_1 \cdot \hat{p}_1 + s_2 \cdot \hat{p}_2\nonumber \\
    & \quad -s_3 \cdot (V-\gamma) \in \Sigma[(e,\hat{x},\hat{u})], \label{eq:A4} \\
    &\overline{u}_i - \kappa_i + s_4^j \cdot (V - \gamma) + s_5^j \cdot \hat{p}_1 \nonumber \\
    &\quad  + s_6^j \cdot \hat{p}_2 \in \Sigma[(e,\hat{x},\hat{u})], \forall j = 1,...,m, \label{eq:A5} \\
    &\kappa_i - \underline{u}_i+ s_7^j \cdot (V - \gamma) + s_8^j \cdot \hat{p}_1 \nonumber \\
    &\quad  + s_9^j \cdot \hat{p}_2 \in \Sigma[(e,\hat{x},\hat{u})], \forall j = 1,...,m. \label{eq:A6}
\end{align}
\endgroup
In the above formulation, SOS polynomials $s_i$ and $s_k^j$ are multipliers used in the generalized S-procedure, and $L_1 = \epsilon_1e'e, L_2 = \epsilon_2 [e;\hat{x};\hat{u}]'[e;\hat{x};\hat{u}]$ with $\epsilon_1$ and $\epsilon_2$ small positive number on the order of $10^{-6}$. Note that $\textbf{sosopt}$ is a nonconvex problem, since the decision variables $\frac{\partial V}{\partial x}$ and $\kappa$, $s$ and $V$ multiply with each other. To tackle this problem, we decompose it into two tractable subproblems to iteratively search between the simulation function $V$ and multipliers and the control law $s_i, s_k^j, \kappa$.

\begin{algorithm}[H] 
\caption{Iterative method}  
\begin{algorithmic}[1]
\Require{A simulation function $V^0$ such that (\ref{eq:A3}) to (\ref{eq:A6}) are feasible for proper choice of $s_i$, $s_k^j$ $\gamma$, $\kappa$.}
\Ensure{($\kappa, \gamma, V$) such that with the volume of $\Omega_{\gamma}^V$ has been shrunk.}
\State $\gamma$-step: minimization problem with $V^0$ fixed,
\begin{align}
    &\min_{\gamma, \kappa, s_1, s_2, s_3, s_4^j,s_5^j,s_6^j,s_7^j,s_8^j,s_9^j} \gamma \nonumber \\
    \textnormal{s.t.} \ & \kappa \in \mathbb{R}^m[(e, \hat{x}, \hat{u})], \ s_i \in \Sigma[(e,\hat{x},\hat{u})], \forall i = 1,...,3, \nonumber \\
    &s_k^j \in \Sigma[(e,\hat{x},\hat{u})], \ \forall j = 1,...,m, \ \forall k = 4,...,9, \nonumber \\
    &-\frac{\partial V^0}{\partial e}\cdot \left(f_e + g_e \cdot \kappa \right)  - L_2 + s_1 \cdot \hat{p}_1 +s_2 \cdot \hat{p}_2 \nonumber \\
    &\quad   - s_3 \cdot (V^0-\gamma) \in \Sigma[(e,\hat{x},\hat{u})], \nonumber \\
    &\overline{u}_i - \kappa_i + s_4^j \cdot (V^0 - \gamma)  + s_5^j \cdot \hat{p}_1 \nonumber \\
    &\quad + s_6^j \cdot \hat{p}_2(\hat{u}) \in \Sigma[(e,\hat{x},\hat{u})], \forall j = 1,...,m, \nonumber \\
    &\kappa_i - \underline{u}_i+ s_7^j \cdot (V^0 - \gamma) + s_8^j \cdot \hat{p}_1 \nonumber \\
    & \quad + s_9^j \cdot \hat{p}_2
    \in \Sigma[(e,\hat{x},\hat{u})], \forall j = 1,...,m.  \nonumber 
\end{align}
\State $V$-step: feasibility problem over decision variables $V, s_1, s_2, s_5^j, s_6^j, s_8^j, s_9^j, s_{10}$:
\begin{align}
    & (s_{10} - \epsilon_3),\ V-L_1 \in \Sigma[e], \ s_1, s_2 \in \Sigma[(e,\hat{x},\hat{u})],   \nonumber \\
    &s_k^j \in \Sigma[(e,\hat{x},\hat{u})], \ \forall j = 1,...,m, \ \forall k = 5,6,8,9, \nonumber \\
    &-s_{10}\cdot (V^0 - \gamma^*) + (V - \gamma^*) \in \Sigma[e], \label{eq:shapefcn} \\
    &-\frac{\partial V}{\partial e} \cdot (f_e + g_e \cdot \kappa^0)  - L_2 + s_1 \cdot \hat{p}_1 + s_2 \cdot \hat{p}_2 \nonumber\\
    &\quad  -s_3^0 \cdot (V-\gamma^*) \in \Sigma[(e,\hat{x},\hat{u})], \nonumber \\
    &\overline{u}_i - \kappa_i^0 + s_4^{j,0} \cdot (V - \gamma^*) + s_5^{j} \cdot \hat{p}_1\nonumber \\
    &\quad + s_6^{j} \cdot \hat{p}_2 \in \Sigma[(e,\hat{x},\hat{u})], \forall j = 1,...,m, \nonumber\\
    &\kappa_i^0 - \underline{u}_i+ s_7^{j,0} \cdot (V - \gamma^*)   + s_8^{j} \cdot \hat{p}_1 \nonumber \\
    &\quad  + s_9^{j} \cdot \hat{p}_2 \in \Sigma[(e,\hat{x},\hat{u})], \forall j = 1,...,m. \nonumber 
\end{align}
\end{algorithmic}
\end{algorithm}

The $\gamma$-step fixes the simulation function $V$ to be $V^0$ computed from the $V$-step of the previous iteration, and tries to find the smallest $\Omega_{\gamma}^{V^0}$ by minimizing $\gamma$ over the control law and multipliers. Since $s_3$ and $\gamma$ multiply with each other, the minimization problem is a bilinear problem and can be effectively solved by bisecting $\gamma$.

In the optimization problem of the $V$-step, $\gamma^*, \kappa^0, s_3^0, s_4^{j,0}$ and $s_7^{j,0}$ are obtained from the $\gamma$-step, and $\epsilon_3$ the small positive number on the order of $10^{-6}$ ensures that $s_{10}$ can't take the value of zero. The constraint (\ref{eq:shapefcn}) enforces the \He{error bound $\Omega_{\gamma^*}^V$ certified by} the $V$-step of current iteration to be contained in $\Omega_{\gamma^*}^{V^0}$. After the $\gamma$-step, the constraints of the $\gamma$-step are active for $V^0$. In the $V$-step, a new feasible $V$ is computed, which is the analytic center of linear matrix inequality constraints. Thus the feasibility problem in the $V$-step pushes $V$ away from the constraints, and gives the next $\gamma$-step more freedom to decrease $\gamma$.

The input to Algorithm 1 is a feasible initial guess for $V^0$. One candidate might be a Lyapunov function $\overline{V}$ obtained by solving Lyapunov equations using the linearized error dynamics with LQR controllers. However, $\overline{V}$ may be too coarse to be feasible for the constraints (\ref{eq:A3}) to (\ref{eq:A6}), which are required to hold on a large set of $\hat{\mathcal{X}}$ and $\hat{\mathcal{U}}$. Here, we propose an algorithm for initializing $V^0$, the idea of which is to enforce the constraints on small subsets of $\hat{\mathcal{X}}, \hat{\mathcal{U}}$, at the beginning and gradually expand the subsets while maintaining feasibility by bisearching between $V$ and $\kappa$.

\begin{algorithm}[H] 
\caption{$V^0$ initialization}  
\begin{algorithmic}[1]
\Require{A Lyapunov function $\overline{V}$ for the linearized closed-loop error dynamics, $\hat{\mathcal{X}}_1 \subset \hat{\mathcal{X}}$ and  $\hat{\mathcal{U}}_1 \subset \hat{\mathcal{U}}$.}

\Ensure{A simulation function $V^0$ such that (\ref{eq:A3}) to (\ref{eq:A6}) are feasible for proper choice of $s_i$, $s_k^j$ $\gamma$, $\kappa$.}
\State $t \gets 1$.
\State $V_1 \gets \overline{V}$.
\While{$\hat{\mathcal{X}}_t \subset \hat{\mathcal{X}}$ \He{or} $\hat{\mathcal{U}}_t \subset \hat{\mathcal{U}}$}
\State Solve for $\kappa_t$ by enforcing (\ref{eq:A3}) to (\ref{eq:A6}) to hold

$\forall \hat{x} \in \hat{\mathcal{X}}_t, \ \forall \hat{u} \in \hat{\mathcal{U}}_t$, while fixing $V = V_t$.
\State Solve for \He{$V_{t+1}$} by enforcing (\ref{eq:A3}) to (\ref{eq:A6}) to hold  

$\forall \hat{x} \in \hat{\mathcal{X}}_t, \ \forall \hat{u} \in \hat{\mathcal{U}}_t$, while fixing $\kappa = \kappa_t$.
\State $\hat{\mathcal{X}}_{t+1} \gets \hat{\mathcal{X}}_{t} \cup \Delta{\hat{\mathcal{X}}}_t$, where $\Delta{\hat{\mathcal{X}}}_t \subset \mathbb{R}^{\hat{n}}$.
\State $\hat{\mathcal{U}}_{t+1} \gets \hat{\mathcal{U}}_{t} \cup \Delta{\hat{\mathcal{U}}}_t$, where $\Delta{\hat{\mathcal{U}}}_t \subset \mathbb{R}^{\hat{m}}$.
\State $t \gets t+1$.
\EndWhile
\State $V^0 \gets V_{N}$.
\end{algorithmic}
\end{algorithm}

\subsection{Outline of the Design Procedure} \label{outline}
We now outline our overall control design process. In particular, if the error bound obtained from the procedure in Section \ref{SOSopt} is too large, further modification to the controller for the abstraction will be necessary. We propose using the following heuristic approach to converge on an allowable error bound:

$(i)$ Select an abstraction $\hat{\Sigma}$ and map $\pi(\cdot)$ based on the application. Then, design a control $\hat{u}$ for $\hat{\Sigma}$ which constrains $\hat{x}(t)$ and $\hat{u}(t)$ to the respective sets $\hat{\mathcal{X}}$ and $\hat{\mathcal{U}}$ (using, e.g., an MPC-based approach).

$(ii)$ Complete the procedure outlined in Section \ref{SOSopt}, resulting in a low level controller $\kappa(e(t), \hat{x}(t), \hat{u}(t))$ for $\Sigma$ and a bound on the error $\Omega_{\gamma}^V$.

$(iii)$ If step $(ii)$ fails to be feasible, or if the error bound is too large, return to step $(i)$ and tighten the constraints on $\hat{x}(t)$ and $\hat{u}(t)$ by shrinking the sets $\hat{\mathcal{X}}$ and $\hat{\mathcal{U}}$. Repeat the process until the error bound is sufficiently small.

To illustrate this procedure, we provide two examples in Sections \ref{example1} and \ref{example2} in which a suitable error bound was obtained by iterating through steps $(i)$ - $(iii)$. In both examples, only a few iterations were necessary.

\section{Vehicle Platoon Example} \label{example1}
\subsection{Model}

\begin{figure}[H]
\centering
\includegraphics[width=0.5\textwidth]{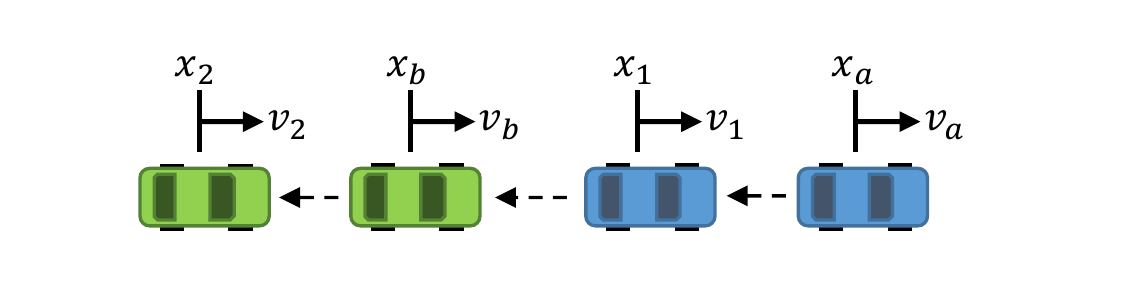}
\vspace{-7.5mm}
\caption{Two platoons of vehicles with leader vehicles $a$ and $b$ and follower vehicles $1$ and $2$. The position and velocity of each vehicle is labelled.}
\label{platoonFigure}
\end{figure}

In this section we apply the previous approach to a vehicle platoon application adapted from \cite{smith2018hierarchical}. We consider the longitudinal dynamics of two vehicle platoons, where the preceding platoon is comprised of a leader vehicle and $n$ follower vehicles, and the latter of one leader vehicle and $m$ follower vehicles (see Figure \ref{platoonFigure} where $n = m = 1$). The task is to formally synthesize a centralized controller for an abstraction of the platoon, with the goal of achieving a high level objective expressed in STL. A low level controller as in Section \ref{sec_overview} is then designed to minimize the error bound between the platoon and its abstraction.

The longitudinal dynamics of the leader vehicle are 
\begin{equation}
\dot{x}_a = v_a, \quad \dot{v}_a = -v^0 + v_a - \rho \cdot (v_a)^2/M + u_a,
\end{equation}
where $x_a$ and $v_a$ are position and velocity, respectively, $v^0$ is a nominal velocity, $\rho$ is an air drag coefficient, $M$ is the vehicle mass, and $u_a$ is the control input to be designed. To improve readability, the time arguments for the variables of differential equation models are omitted in the examples. Furthermore, we note that the quadratic air drag term was not included in \cite{smith2018hierarchical}, which specialized to affine systems. The follower vehicles are governed by a simple vehicle-following model
\begin{align}
\dot{x}_1 &= v_1, \nonumber \\
\dot{v}_1 &= v^0 - v_1 - \rho \cdot (v_1)^2/M + g(x_a - x_1), \nonumber \\
\dot{x}_i &= v_i, \nonumber \\
\dot{v}_i &= v^0 - v_i - \rho \cdot (v_i)^2/M + g(x_{i-1} - x_i),
\end{align}
for $i = 2, \dots, n$, where $g: \R \to \R$ is a monotone increasing function. The function $g(\cdot)$ models a unilateral spring force acting on vehicle $i$. Its purpose is to preserve an inter-vehicle spacing of $\ell \in \R$ within the platoon, and thus $g(\ell) = 0$.

The second vehicle platoon has the same dynamics as the first; however, as a safety measure we also model a spring force between the $n$th follower vehicle and platoon leader $b$:
\begin{align} \label{withSpring}
\dot{x}_b &= v_b, \nonumber \\
\dot{v}_b &= -v_b + v^0 - \rho \cdot (v_b)^2/M + g(x_n - x_b) + u_b.
\end{align}

We note that as long as $g(\cdot)$ is polynomial, the platoon dynamics are polynomial, and are therefore amenable to SOS programming. For simplicity, we use a linear spring force
\begin{equation} \label{springForce}
g(s) = k \cdot (s - \ell),
\end{equation}
where $k \in \R$ is a spring constant. We note that the overall platoon model may be represented compactly as in \eqref{concrete}.

\subsection{Abstraction}
To obtain an abstract model of the platoon, we consider only the dynamics of the two leader vehicles
\begin{align}
\dot{\hat{x}}_a &= v_a, \quad \dot{\hat{v}}_a = -\hat{v}^0 + \hat{v}_a - \rho \cdot (\hat{v}_a)^2/M + \hat{u}_a, \\
\dot{\hat{x}}_b &= \hat{v}_b, \quad \dot{\hat{v}}_b = -\hat{v}^0 + \hat{v}_b - \rho \cdot (\hat{v}_b)^2/M + \hat{u}_b,
\end{align}
for which we will formally synthesize a control strategy. We note, in particular, that this model ignores the spring force in \eqref{withSpring}, as the position of vehicle $n$ is not accessible. This effectively introduces a disturbance into the error dynamics \eqref{errDyn}, which is compensated for by the interface \eqref{interface} - for further details see \cite{smith2018hierarchical}. Lastly, we define the manifold \eqref{manifold} as follows
\begin{alignat*}{2}
x_a & = \hat{x}_a, \quad x_i = \hat{x}_a - i\ell, \ i = 1, \dots, n, \\
v_a & = v_1 = \dots = v_n = \hat{v}_a, \\
x_b & = \hat{x}_b, \quad x_{n+i} = \hat{x}_b - i\ell, \ i = 1, \dots, m, \\
v_b & = v_{n+1} = \dots = v_{n+m} = \hat{v}_b,
\label{aggM}
\end{alignat*}
which is affine, and therefore amenable to our approach. Physically, this manifold means that 1) the leader vehicles match their respective reference positions, with all follower vehicles leaving $\ell$ meters of distance in front of them, and 2) all vehicles within each platoon match their respective reference velocity.

\begin{figure*}
    \centering
    \input{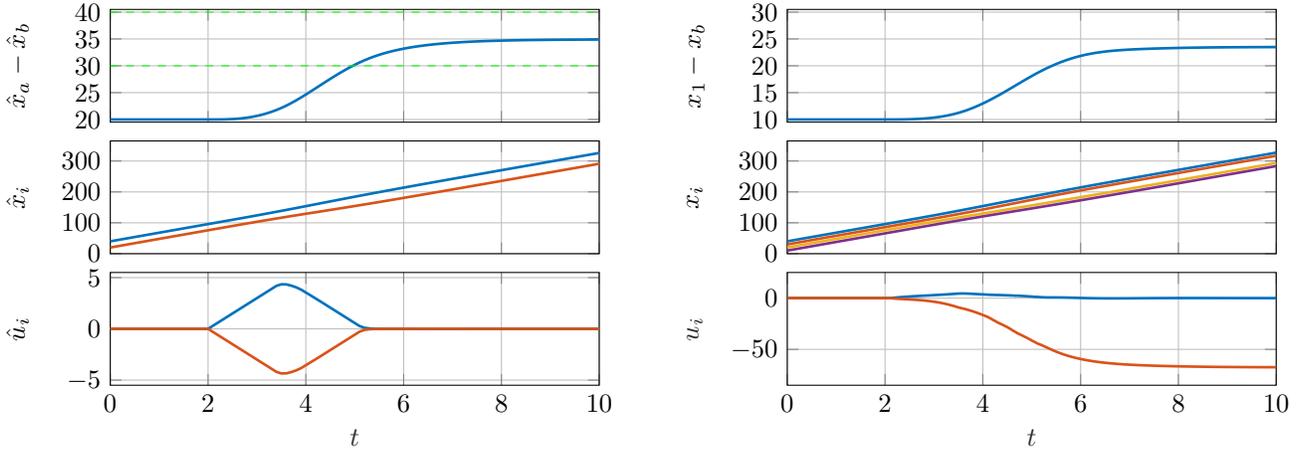}
    \caption{Simulation results for the platoon separation example. In the abstraction (left), the leader vehicles increase their distance to reach the desired set, indicated by the green dashed lines. The gap maneuver is successfully realized on the full platoon (right), resulting in an increased distance between platoon leader $b$ and follower vehicle $1$.}
    \label{platoon_result}
\end{figure*}

\subsection{Specification: Platoon Separation}
We consider a two platoons of vehicles each with one leader and one follower vehicle ($n = m = 1$). A common control objective for platooning is a gap widening maneuver, with the goal of allowing another vehicle to merge into the lane occupied by the platoon. This objective may be expressed in STL as follows
\begin{equation} \label{platoonSTL}
\square (\varphi \wedge (S \implies \lozenge_{[0, T]} \psi )).
\end{equation}
Here, the $\square$ and $\lozenge$ symbols are the "always" and "eventually" operators; thus, $\varphi$ expresses state constraints which must be satisfied at all times
\begin{align} \label{safetyConstr}
& \varphi := (|\hat{v}_a - 28| \leq 5) \wedge (|\hat{v}_b - 28| \leq 5) \wedge (|\hat{u}_a| \leq 8) \nonumber \\
& \qquad \wedge (|\hat{u}_b| \leq 8) \wedge \ (|\hat{x}_a - \hat{x}_b - 20| \ \leq 20),
\end{align}
and $\psi$ specifies a desired gap between the leader vehicles
\begin{align} \label{reachSet}
\psi := (|\hat{x}_a - \hat{x}_b - 35| \leq 5),
\end{align}
which must be attained within $T$ time-steps after the binary signal $S$ becomes equal to $1$. In practice, $S$ could represent another vehicle indicating that it would like to merge in between follower vehicle $n$ and leader vehicle $b$ in the platoon. Using the MPC-based approach in \cite{raman2017}, we designed a controller for a time-discretized model of the abstraction which achieves \eqref{platoonSTL}.

The constraints in \eqref{safetyConstr} are easily converted to semi-algebraic constraint sets $\hat{\mathcal{X}}$ and $\hat{\mathcal{U}}$. Decision variables $V$ and $\kappa$ are characterized by degree-2 polynomials with unknown parameters. In this example, $e(0)$ is assumed to be chosen within $\Omega_{\gamma}^V$, therefore we can relax (\ref{eq:A2}) to be
\begin{align}
    \dot{V}(e, \hat{x}, \hat{u}) < 0, \ \forall e, \hat{x}, \hat{u}, \ \textnormal{s.t.} \ V(e) = \gamma, \ \hat{x} \in \hat{\mathcal{X}}, \ \hat{u} \in \hat{\mathcal{U}}. \nonumber
\end{align}
In this case, $V$ becomes a barrier function, and $\Omega_{\gamma}^V$ is a forward invariant set for $e(t)$. Consequently, we have a weaker constraint for $s_3$ than (\ref{eq:A3}): $s_3(e,\hat{x},\hat{u}) \in \R[(e, \hat{x}, \hat{u})]$. The SOS programming was formulated and translated into SDP using the sum-of-square module in SOSOPT \cite{Pete:13} on Matlab, and solved by Mosek \cite{Mosek:17}. The computation was performed on a workstation with a 2.7 [GHz] Intel Core i5 64 bit processor and 8 [GB] of RAM, and this 8-state example takes around 15 minutes to solve. The simulation results are shown and explained in Figure \ref{platoon_result}.

\section{Double Pendulum Example} \label{example2}
\begin{figure}[h]
	\centering
	\includegraphics[width=0.24\textwidth]{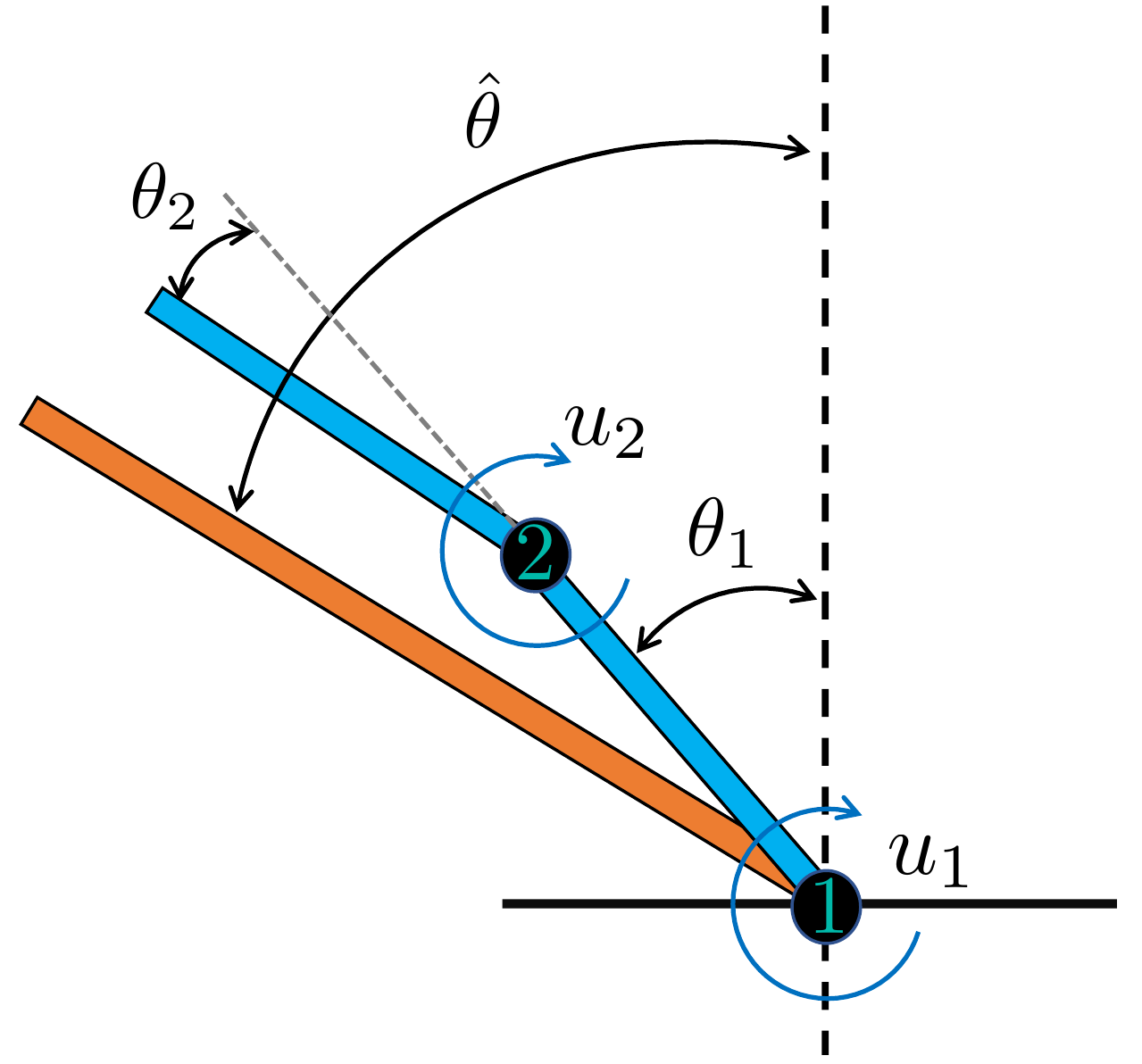}
	\caption{Double pendulum and its abstraction as a single pendulum. The angular positions are labelled.}
	\label{fig:angles}    
\end{figure}

In this example we consider the fully-actuated double pendulum dynamics, where torques are applied at the both joints. The polynomial dynamics for the double pendulum shown below are obtained by a least-squares approximation of the full equation for $(x_1, x_3) \in [-1, 1] \times [-1, 1]$:
\begin{align}
\bmat{\dot{x}_1 \\ \dot{x}_2 \\ \dot{x}_3 \\ \dot{x}_4} &= \bmat{x_2 \\ f_2(x_1, x_2, x_3, x_4) \\ x_4 \\ f_4(x_1, x_2, x_3, x_4)} + \bmat{0 & 0 \\ 8 & -31.2 \\ 0 & 0 \\ -31.2 & 391.2}\bmat{u_1 \\ u_2}, \nonumber
\end{align}
\begingroup
\begin{align}
f_2 &=   - 3.447 x_1^3 + 2.350 x_1^2 x_3 + 1.303 x_1 x_3^2
  + 3.939 x_3^3 \nonumber \\
  & \quad \quad + 21.520 x_1 - 5.000 x_3, \nonumber \\
f_4 &=  4.023 x_1^3 - 36.551 x_1^2 x_3 - 4.131 x_2^2 x_3 - 27.060 x_3^3 \nonumber \\
& \quad \quad - 25.115 x_1 + 77.700 x_3, \nonumber
\end{align}
\endgroup
where $x_1$ and $x_3$ represent $\theta_1$ and $\theta_2$ (shown in Fig. \ref{fig:angles}), which are angular positions of the first and second links (relative to the first link); $x_2$ and
$x_4$ represent $\omega_1$ and $\omega_2$, which are angular velocities of two links; $u_1$ and $u_2$ are torques applied at the joint 1 and joint 2 (shown in Fig. \ref{fig:angles}), respectively. 

The control objective is to bring $[\theta_1; \omega_1]$ from $[-0.5; 0.5]$ to $[0.3; 0]$ and to maintain it there, while minimizing the relative angular position $\theta_2$ and relative angular velocity $\omega_2$ between the two links. We choose the single inverted pendulum model as the abstraction to the original dynamics, design controllers for the abstraction and then synthesize the low level controller $\kappa$ to minimize the error bound between the concrete and abstract models. The polynomial dynamics of the abstraction are computed by a least-squares approximation of the full equation for $\hat{x}_1 \in [-1, 1]$:
\begin{align}
    \bmat{\dot{\hat{x}}_1 \\ \dot{\hat{x}}_2} = \bmat{\hat{x}_2 \\ -5.131 \hat{x}_1^3 + 32.1 \hat{x}_1} + \bmat{0 \\ 9.1}\hat{u}, \nonumber
\end{align}
where $\hat{x}_1$ represents $\hat{\theta}$ (shown in Fig. \ref{fig:angles}), the angular position of the single pendulum, $\hat{x}_2$ represents $\hat{\omega}$, the angular velocity, and $\hat{u}$ represents the torque applied at the joint~1. Based on the control objective, the manifold is chosen as:
\begin{align}
    x = P \hat{x}, \ \text{where} \ P = \bmat{\boldsymbol{I}^2; \mathbf{0}^{2 \times 2}}. \nonumber
\end{align}

\begin{figure}[h]
	\centering
	\begin{subfigure}[b]{0.45\textwidth}
		\centering
		\includegraphics[width=1.1\textwidth]{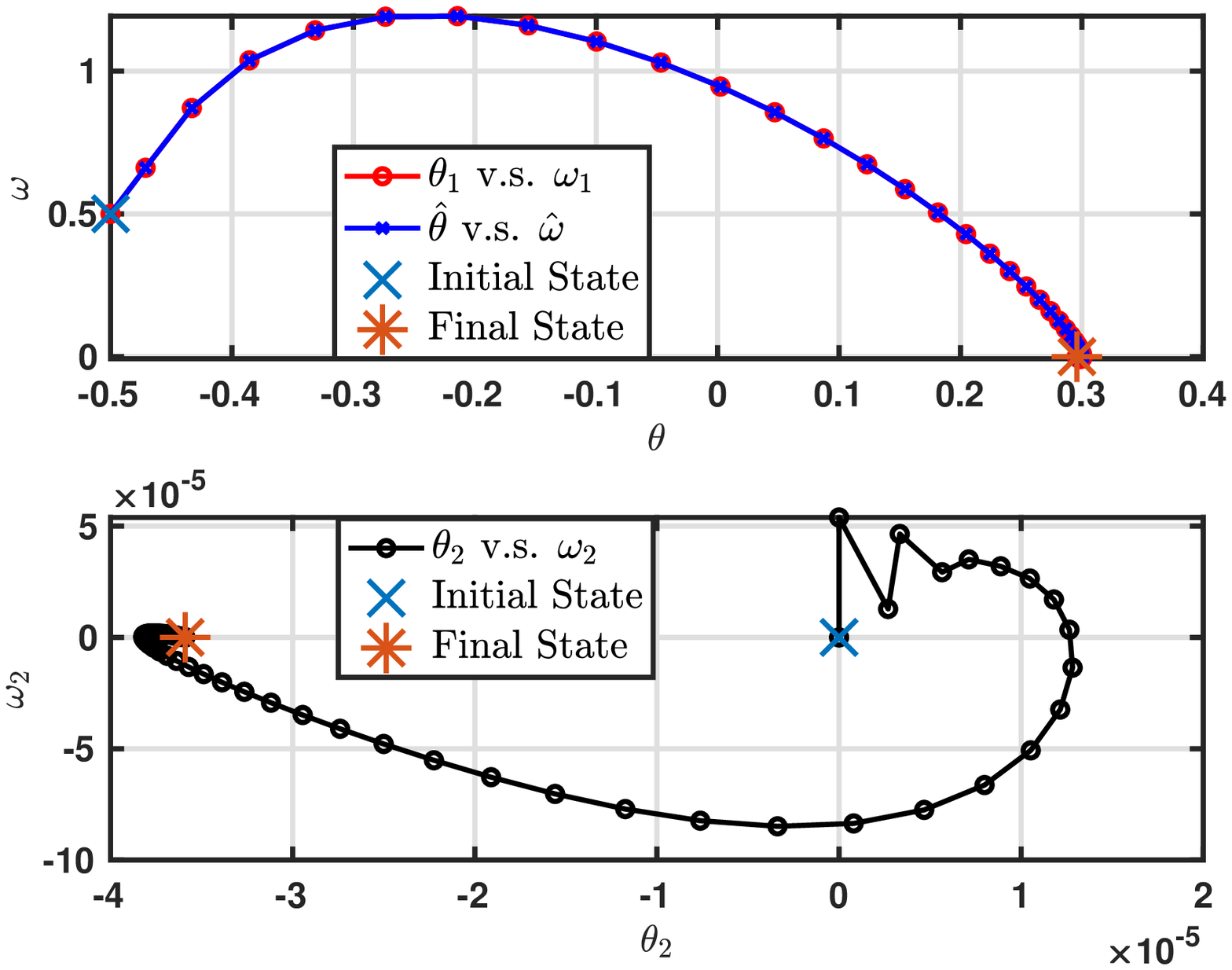}
		\caption{Phase portraits for the double pendulum example}  
	\end{subfigure}
	\hfill
	\begin{subfigure}[b]{0.45\textwidth}  
		\centering 
		\includegraphics[width=1.1\textwidth]{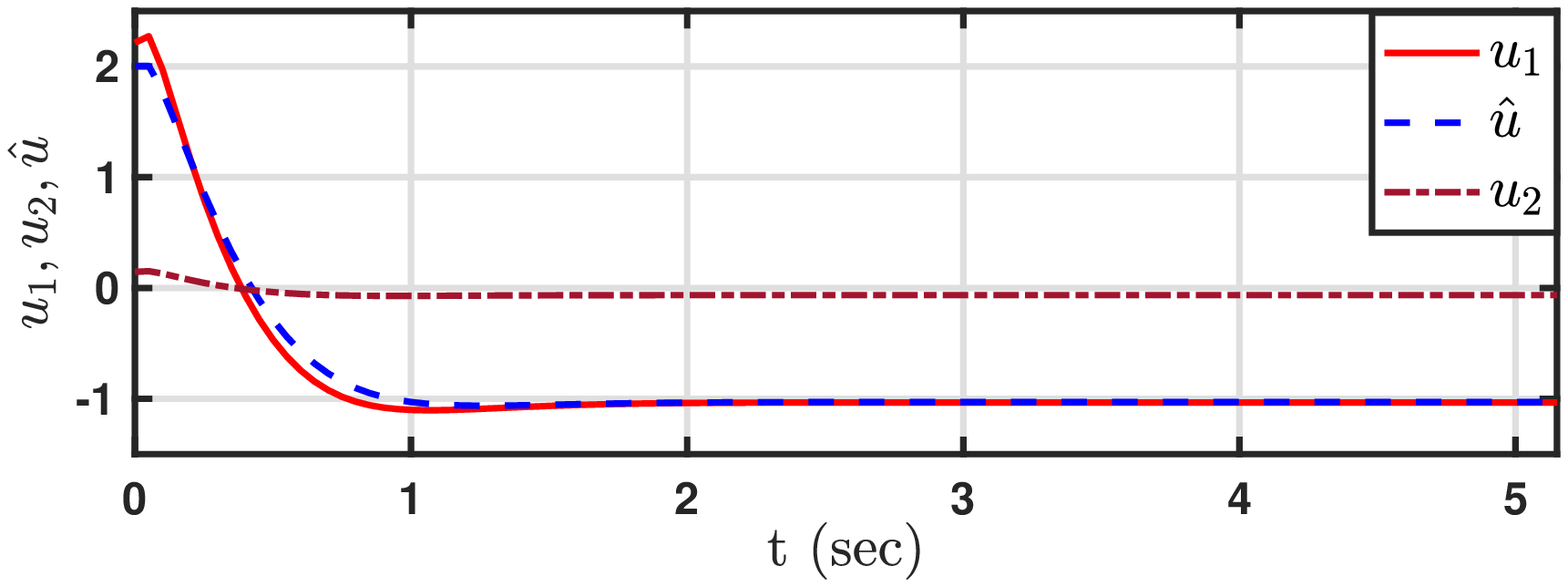}
		\caption{Control inputs vs. time}
	\end{subfigure}
	\caption{Simulation results for the double pendulum example}
	\label{fig:Pendulum}    
\end{figure}

To synthesize the low level control law $\kappa$ and \He{compute} the error bound, we enforce the following states and inputs constraints when design the MPC controller for the abstract system
\begin{align}
    \hat{x}(t) \in \hat{\mathcal{X}} &= \{\hat{x} \in \mathbb{R}^{\hat{n}} : \vert \hat{x}_1 \vert \leq 0.6, \ \vert\hat{x}_2 \vert \leq 1.3\}, \nonumber \\
    \hat{u}(t) \in \hat{\mathcal{U}} &= \{\hat{u} \in \mathbb{R}^{\hat{m}} : \vert \hat{u} \vert \leq 5\}. \nonumber
\end{align} 
$V$ and $\kappa$ are parameterized by degree-2 and degree-4 polynomials, respectively. The whole computation including the initialization takes less than 5 minutes to perform on a workstation  with a 2.7 [GHz] Intel Core i5 64 bit processor and 8 [GB] of RAM. The simulation results are shown in Figure \ref{fig:Pendulum}. We can see that the error between the phase portrait of $\theta_1$ and $\omega_1$ and that of $\hat{\theta}$ and $\hat{\omega}$ is small, the phase portrait of $\theta_2$ and $\omega_2$ stays closely to the origin, and the simulation trajectory of $u_1$ is close to that of $\hat{u}$.


% ...



\section{Conclusion} \label{sec_conclusion}
In this paper we proposed an abstraction-based control design approach for nonlinear control-affine systems. In particular, we utilized SOS optimization to compute a low-level controller and a bound on error states between the original system and its abstraction. This optimization procedure assumes that the dynamics are approximated with a polynomial model and the abstract state and input are restricted to semi-algebraic constraint sets. We also proposed an iterative procedure to refine these constraint sets to the necessary level of restrictiveness so that an acceptable error bound is found. Finally, we applied our methods to two examples which demonstrate its broad applicability.





\bibliographystyle{IEEEtran}
\bibliography{IEEEabrv,bibfile}

\end{document}